\documentclass[12pt]{amsproc}
\usepackage{amsmath,amssymb,amscd}

\begin{document}

\title[$K$-Theory and Perturbations]{$K$-Theory and Perturbations of Absolutely Continuous Spectra}
\author{Dan-Virgil Voiculescu}
\address{D.V. Voiculescu \\ Department of Mathematics \\ University of California at Berkeley \\ Berkeley, CA\ \ 94720-3840}
\thanks{Research supported in part by NSF Grant DMS-1301727.}
\date{}

\begin{abstract}
We study the $K_0$-group of the commutant modulo a normed ideal of an $n$-tuple of commuting Hermitian operators in some of the simplest cases. In case $n=1$, the results, under some technical conditions are rather complete and show the key role of the absolutely continuous part when the ideal is the trace-class. For a commuting $n$-tuple, $n \ge 3$ and the Lorentz $(n,1)$ ideal, we show under an absolute continuity assumption that the commutant determines a canonical direct summand in $K_0$. Also, certain properties involving the compact ideal, established assuming quasicentral approximate units mod the normed ideal, have weaker versions which hold assuming only finiteness of the obstruction to quasicentral approximate units.
\end{abstract}

\maketitle

\section{Introduction}
\label{sec1}

By the Kato--Rosenblum theorem (see \cite{8} or \cite{9}) the absolutely continuous part is preserved up to unitary equivalence under trace-class perturbations. We generalized this in \cite{11} to perturbations of commuting $n$-tuples of Hermitian operator, under the requirements that the commutation be preserved and that the normed ideal be the $(n,1)$-Lorentz ideal, using an approach based on the realization that such questions ultimately boil down to existence or non-existence of quasicentral approximate units relative to the ideal (for more see also \cite{1}, \cite{12}, \cite{13}). Recently we found that an object which appears naturally in these perturbation questions is the commutant modulo a normed ideal of an $n$-tuple of operators (\cite{3}, \cite{14}, \cite{15}, \cite{16}). In this vein, in this note we show that results about perturbations of absolutely continuous spectra have consequences for the $K$-theory of the corresponding commutants mod normed ideals. Note also applications of our approach to absolutely continuous spectra (\cite{11}) in non-commutative geometry in \cite{5}.

The paper has five sections including this introduction. Section~\ref{sec2} deals with preliminaries and in particular emphasizes how equivalences of projections in commutants mod normed ideals amount to unitary equivalence mod the normed ideal of the compressions of the defining operators. Section~\ref{sec3} deals with the case of one hermitian operator and the perturbation theory provides a rather complete picture of what to expect from the $K_0$ group of the commutant mod an arbitrary ideal. In essence in the case of the trace-class the multiplicity function of an absolutely continuous part describes the whole group, while for other ideals, under suitable restrictions $K_0$ is trivial. In Section~\ref{sec4} we show that for commuting $n$-tuples of Hermitian operators, $n \ge 3$ with Lebesgue absolutely continuous spectral measure, the $K_0$-group of the von~Neumann algebra which is the commutant of the $n$-tuple appears as a direct summand in the $K_0$-group of the commutant mod the Lorentz-type $(n,1)$-ideal (\cite{7}, \cite{10}) in a canonical way (that is independent of replacing the $n$-tuple by an equivalent mod the ideal one). One should, however, be aware that for $n > 1$, there may be other $K_0$-classes like ones related to cyclic cohomology. Section~\ref{sec5} is about technical results involving the compact ideal assuming only finiteness of the obstruction instead of existence of quasicentral approximate units relative to the ideal. In particular these results apply quite generally when the ideal is the Macaev ideal.

\section{Preliminaries}
\label{sec2}

\noindent
{\bf 2.1.} With ${\mathcal H}$ denoting throughout this paper a separable infinite-dimensional Hilbert space, if $\tau = (T_j)_{1 \le j \le n}$ is an $n$-tuple of bounded Hermitian operators and $({\mathcal J},|\ |_{\mathcal J})$ is a normed ideal of compact operators, we shall consider the commutant $\mod {\mathcal J}$ of $\tau$, which is the $*$-algebra (see \cite{15}, \cite{16})
\[
{\mathcal E}(\tau;{\mathcal J}) = \{X \in {\mathcal B}({\mathcal H}) \mid [X,T_j] \in {\mathcal J},\ 1 \le j \le n\}
\]
with the norm
\[
\||X\|| = \|X\| + \max_{1 \le j \le n} |[T_j,X]|_{\mathcal J}.
\]
If $\tau' = (T'_j)_{1 \le j \le n}$ is so that $T'_j - T_j \in {\mathcal J}$, $1 \le j \le n$, then ${\mathcal E}(\tau;{\mathcal J}) = {\mathcal E}(\tau';{\mathcal J})$ and the norms $\||\ \||$ and $\||\ \||'$ are equivalent, actually $\||X\|| \le \||X\||' + C\|X\|$ and $\||X\||' \le \||X\|| + C\|X\|$ for some $C > 0$. With ${\mathcal K}$ denoting the compact operators on ${\mathcal H}$, we use the notation ${\mathcal K}(\tau;{\mathcal J}) = {\mathcal E}(\tau;{\mathcal J}) \cap {\mathcal K}$ and ${\mathcal E}/{\mathcal K}(\tau;{\mathcal J}) = {\mathcal E}(\tau;{\mathcal J})/{\mathcal K}(\tau;{\mathcal J})$. Clearly in a purely algebraic way there is a homomorphism $p: {\mathcal E}/{\mathcal K}(\tau;{\mathcal J}) \to {\mathcal B}/{\mathcal K}$ identifying ${\mathcal E}/{\mathcal K}(\tau;{\mathcal J})$ with a $*$-subalgebra of the Calkin algebra ${\mathcal B}/{\mathcal K}$.

\bigskip
\noindent
{\bf 2.2.} Let ${\mathcal M}_n$ denote the $n \times n$ matrices. Then we have an obvious identification
\[
{\mathcal M}_n{\mathcal E}(\tau;{\mathcal J}) \sim {\mathcal E}(\tau \otimes {\mathcal I}_n;{\mathcal M}_n{\mathcal J}).
\]
Here $\tau \otimes {\mathcal I}_n$ is viewed as acting on ${\mathcal H} \otimes {\mathbb C}^n \simeq {\mathcal H}^n$ and ${\mathcal M}_n {\mathcal J}$ is endowed with one of the several equivalent natural norms. Perhaps the best choice of norm on ${\mathcal M}_n{\mathcal J}$ is to use an isometry $V: {\mathcal H} \to {\mathcal H}^n$ and consider $|V^* \cdot V|_{\mathcal J}$, though choosing another equivalent norm will be irrelevant in what follows.

\bigskip
\noindent
{\bf 2.3.} We will also use the invariant
\[
k_{\mathcal J}(\tau) = \liminf_{A \in {\mathcal R}^+_1({\mathcal H})} \max_{1 \le j \le n} |[A,T_j]|
\]
where ${\mathcal R}^+_1({\mathcal H})$ denotes the finite rank positive contractions on $A$ endowed with its natural order, with respect to which the liminf is taken (see \cite{11}, \cite{12}, \cite{13}). If the $n$-tuples $\tau$ and $\tau'$ are equal $\text{mod } {\mathcal J}$ and the finite rank operators ${\mathcal R}$ are dense in ${\mathcal J}$, then clearly $k_{\mathcal J}(\tau) = k_{\mathcal J}(\tau')$.

\bigskip
\noindent
{\bf 2.4.} We shall denote by $({\mathcal C}_p,|\ |_p)$ the Schatten--von~Neumann $p$-classes, $1 \le p < \infty$ (see \cite{7}, \cite{10}). The $(p,1)$-Lorentz type ideal will be denoted by $({\mathcal C}_p^-,|\ |_p^-)$ where $|X|^-_p = \sum_{j \ge 1} s_j(X)j^{-1+1/p}$ with $s_1(X) \ge s_2(X) \ge \dots$ denoting the eigenvalues of $(X^*X)^{1/2}$ in decreasing order (multiple eigenvalues repeated according to multiplicity). Remark that $({\mathcal C}_1,|\ |_1) = ({\mathcal C}^-_1,|\ |^-_j)$. If ${\mathcal J}$ is either ${\mathcal C}_p$ or ${\mathcal C}_p^-$, we denote $k_{\mathcal J}(\tau)$ by $k_p(\tau)$ and respectively by $k^-_p(\tau)$. If $p > 1$, then $k_p(\tau) \in \{0,\infty\}$, while $k^-_p(\tau)$ may also take finite non-zero values (\cite{1}, \cite{11}, \cite{12}, \cite{13}).

\bigskip
\noindent
{\bf 2.5.} The idempotents in ${\mathcal E}(\tau;{\mathcal J})$ will be denoted by ${\mathcal P}{\mathcal E}(\tau;{\mathcal J})$ and the Hermitian idempotents, that is the orthogonal projection operators, by ${\mathcal P}_h{\mathcal E}(\tau;{\mathcal J})$. If $P \in {\mathcal P}{\mathcal E}(\tau;{\mathcal J})$ then the associated orthogonal projection onto $P{\mathcal H}$, that is $P_{P{\mathcal H}} \in {\mathcal P}_h{\mathcal E}(\tau;{\mathcal J})$. This is a consequence of the smooth functional calculus in ${\mathcal E}(\tau;{\mathcal J})$, that is if $X = X^* \in {\mathcal E}(\tau;{\mathcal J})$ and $f: {\mathbb R} \to {\mathbb R}$ is sufficiently smooth (infinitely differentiable for instance) then $f(X) \in {\mathcal E}(\tau;{\mathcal J})$. Indeed $P_{P{\mathcal H}} = f(PP^*)$ for a suitable function $f$.

In ${\mathcal P}{\mathcal E}(\tau;{\mathcal J})$ there is the equivalence relation $P \sim Q$ if there are $X,Y \in {\mathcal E}(\tau;{\mathcal J})$ so that $P = XY$, $Q = YX$. Then $P \sim P_{P{\mathcal H}}$. Moreover, if $P,Q \in {\mathcal P}_h{\mathcal E}(\tau;{\mathcal J})$, then $P \sim Q$ iff there is a partial isometry $V \in {\mathcal E}(\tau;{\mathcal J})$ so that $VV^* = P$, $V^*V = Q$ (again a consequence of the smooth functional calculus). This means that $K_0({\mathcal E}(\tau;{\mathcal J}))$ can be defined like for $C^*$-algebras using orthogonal projections and partial isometries (\cite{2}, \cite{6}).

\bigskip
\noindent
{\bf 2.6.} Let $P,Q \in {\mathcal P}_h{\mathcal E}(\tau;{\mathcal J})$ and $X \in {\mathcal B}({\mathcal H})$ so that $PXQ = X$. Then we have $X \in {\mathcal E}(\tau;{\mathcal J})$ iff $PT_jPX - XQT_jQ \in {\mathcal J}$, $1 \le j \le n$. In particular, we have $P \sim Q$ iff there is a partial isometry $V$ so that $VV^* = P$, $V^*V = Q$ and $VQT_jQ - PT_jPV \in {\mathcal J}$, $1 \le j \le n$. Note also that $P{\mathcal E}(\tau;{\mathcal J})P \mid P{\mathcal H} = {\mathcal E}(P\tau \mid P{\mathcal H};{\mathcal J})$.

\bigskip
\noindent
{\bf 2.7.} If $\tau$ is an $n$-tuple of commuting Hermitian operators we denote by ${\mathcal H}_{ac}(\tau)$ the Lebesgue-absolutely continuous subspace and by $E_{ac}(\tau)$ the orthogonal projection onto ${\mathcal H}_{ac}(\tau)$. The multiplicity function of the absolutely continuous part of $\tau$, that is of $\tau \mid {\mathcal H}_{ac}$ will be denoted by $mac(\tau)$. Here $mac(\tau)$ is an Lebesgue almost everywhere defined measurable function on ${\mathbb R}^n$ taking values in $\{0,1,2,\dots,\infty\}$. We shall denote by $\omega(\tau) = mac(\tau)^{-1}(\infty)$ which is defined up to null sets. We shall also use the notation ${\mathcal H}_{\text{sing}}(\tau)$ and $E_{\text{sing}}(\tau)$Q for ${\mathcal H} \ominus {\mathcal H}_{ac}(\tau)$ and $I - E_{ac}(\tau)$; the singular subspace and the projection onto it.

\bigskip
\noindent
{\bf 2.8.} We also recall that the essential spectrum of an operator or of an $n$-tuple of commuting Hermitian operators is the spectrum of the image of the operator or $n$-tuple in the Calkin algebra. We shall use the notation $\sigma_e(T)$, $\sigma_e(\tau)$.

\section{Results in the case of one Hermitian operator}
\label{sec3}

We shall denote by ${\mathcal F}(T)$, where $T$ is a Hermitian operator, the additive group of equivalence classes up to almost everywhere equality of measurable functions $f: {\mathbb R}\backslash \omega(T) \to {\mathbb Z}$ so that $|f(t)| \le C\ mac(T) \mid {\mathbb R}\backslash \omega(T)$ for some constant $C$ depending on $f$.

\bigskip
\noindent
{\bf Theorem 3.1.} {\em If $T$ is a Hermitian operator on ${\mathcal H}$, then there is a unique homomorphism
\[
fmac(T): K_0({\mathcal E}(T;{\mathcal C}_1) \to {\mathcal F}(T)
\]
so that if $P \in {\mathcal P}_h({\mathcal E}(T \otimes I_n;{\mathcal C}_1))$ then $fmac(T)([P]_0) = mac(P(T \otimes I_n) \mid P{\mathcal H}^n)$. The homomorphism $fmac(T)$ is surjective.
}

\bigskip
\noindent
{\bf {\em Proof.}} Notice first that if $P \sim Q$, where 
\[
P,Q \in {\mathcal P}_n({\mathcal E}(T \otimes I_n,{\mathcal C}_1))
\]
then $mac(P(T \otimes I_n) \mid P{\mathcal H}^n) = mac(Q(T \otimes I_n) \mid Q{\mathcal H}^n)$. Indeed, by 2.6. there is an isometry from $Q{\mathcal H}^n$ to $P{\mathcal H}^n$ which intertwines $Q(T \otimes I_n) \mid Q{\mathcal H}$ and $P(T \otimes I_n) \mid P{\mathcal H}$ up to ${\mathcal C}_1$. It follows then from the Kato--Rosenblum theorem that $mac(P(T \otimes I_n) \mid P{\mathcal H}^n) = mac(Q(T \otimes I_n) \mid Q{\mathcal H}^n)$ almost everywhere.

Next, remark that if $P \in {\mathcal P}_h({\mathcal E}(T \otimes I_n,{\mathcal C}_1))$ we have
\[
mac(P(T \otimes I_n) \mid P{\mathcal H}^n) \le n\ mac(T)
\]
almost everywhere on ${\mathbb R}$. This is a consequence of the Kato--Rosenblum theorem applied to $T \otimes I_n$ and
\[
P(T \otimes I_n)P + (I \otimes I_n - P)(T \otimes I_n)(I \otimes I_n - P)
\]
after observing that the difference of these operators is in ${\mathcal C}_1$ and that we have almost everywhere that
\[
\begin{aligned}
nmac(T) &= mac(T \otimes I_n) = mac(P(T \otimes I_n)P) \\
&+ mac((I \otimes I_n - P)(T \otimes I_n)(I \otimes I_n-P)) \\
&\ge mac(P(T \otimes I_n)P).
\end{aligned}
\]

The last remark means that
\[
mac(P(T \otimes I_n) \mid P{\mathcal H}^n) \mid {\mathbb R}\backslash\omega(T) \in {\mathcal F}(T).
\]
It is also easily seen that if $P_1 \in {\mathcal P}_h({\mathcal E}(T \otimes I_n,\tau_1))$ and $P_2 \in {\mathcal P}_n({\mathcal E}(T \otimes I_m,{\mathcal C}_1))$ then 
\[
P_1 \oplus P_2 \in {\mathcal P}_h({\mathcal E}(T \otimes I_{n+m},{\mathcal C}_1)
\]
and
\[
\begin{aligned}
&mac((P_1 \oplus P_2)(T \otimes I_{n+m})(P_1 \oplus P_2)) \\
&= mac(P_1(T \otimes I_n)P_1) + mac(P_2(T \otimes I_m)P_2).
\end{aligned}
\]

Having made these observations, to conclude the proof of the existence of $fmac(T)$ it will suffice to see what happens if $P,Q \in {\mathcal P}_h({\mathcal E}(T \otimes I_n,{\mathcal C}_1))$ are so that their $K$-theory classes are equal $[P]_0 = [Q]_0$. Since this is equivalent to $P \oplus I_m \sim Q \oplus I_m$ for some $m \in {\mathbb N}$, we get that
\[
mac(P(T \otimes I_n)P) + m\ mac(T) = mac(Q(T \otimes I_n)Q) + m\ mac(T)
\]
almost everywhere. Since $\omega(T)$ is the set where $mac(T)$ is infinite, we see that on the complement ${\mathbb R}\backslash\omega(T)$, we have the almost everywhere equality
\[
mac(P(T \otimes I_n)P) \mid {\mathbb R}\backslash\omega(T) = mac(Q(T \otimes I_n)Q) \mid {\mathbb R}\backslash\omega(T).
\]
Thus the assignment of
\[
mac(P(T \otimes I_n)P) \mid {\mathbb R}\backslash\omega(T)
\]
to $[P]_0$ is well-defined and gives an additive homomorphism from the semigroup of classes of projections in $K_0({\mathcal E}(T;{\mathcal C}_1))$ into ${\mathcal F}(T)$. Clearly this extends then to a homomorphism of the $K_0$-group and the homomorphism is unique $K_0$ being generated by the classes $[P]_0$ as above.

To prove surjectivity of $fmac(T)$, since ${\mathcal F}(T)$ is generated by its positive functions, it will suffice to show that these are in the range of $fmac(T)$. If $g: {\mathbb R} \to {\mathbb Z}$ is a non-negative Lebesgue measurable function so that $|g(t)| \le n\ mac(T)(t)$ a.e.\ for some $n \in {\mathbb N}$, then there is a $T \otimes I_n$ invariant subspace ${\mathcal X}$ of $({\mathcal H}_{ac}(T))^n \simeq {\mathcal H}_{ac}(T) \otimes {\mathbb C}^n \simeq {\mathcal H}_{cc}(T \otimes I_n)$ so that $mac((T \otimes I_n)({\mathcal X}) = g$. Hence if $P$ is the projection of ${\mathcal H}^n$ onto ${\mathcal X}$, then $fmac(T)([P]_0) = g \mid {\mathbb R}\backslash\omega(T)$.\qed

\bigskip
\noindent
{\bf Theorem 3.2.} {\em Let $T$ be a Hermitian operator with $\sigma(T) = \sigma_e(T)$. Then the homomorphism $fmac(T): K_0({\mathcal E}(T;{\mathcal C}_1)) \to {\mathcal F}(T)$ in Theorem~$3.1.$ is an isomorphism.
}

\bigskip
\noindent
{\bf {\em Proof.}} We first prove the following: if $P,Q \in {\mathcal P}_h({\mathcal E}(T \otimes I_n,{\mathcal C}_1))$ are so that $mac(P(T \otimes I_n) \mid P{\mathcal H}^n) = mac(Q(T \otimes I_n) \mid Q{\mathcal H}^n)$, then in ${\mathcal P}_h({\mathcal E}(T \otimes I_{n+1},{\mathcal C}_1))$ we have $P \oplus I \sim Q \oplus I$. We have $P(T \otimes I_n)P + (I \otimes I_n - P)(T \otimes I_n)(I \otimes I_n - P) \in T \otimes I_n + {\mathcal C}_1$ which using \cite{4} implies the existence of $K' \in {\mathcal B}(P{\mathcal H}^n)$, $K' \in {\mathcal C}_1$ so that $\sigma(P(T \otimes I_n) \mid P{\mathcal H}^n + K') \subset \sigma_e(T)$. Then $S_1 = (P(T \otimes I_n) \mid P{\mathcal H}^n + K') \mid {\mathcal H}_{\text{sing}}(P(T \otimes I_n) \mid P{\mathcal H}^n + K')$ has $\sigma(S_1) \subset \sigma_e(T)$ and since its spectral measure is singular $k_1(S_1) = 0$. Using for instance the adapted non-commutative Weyl--von~Neumann type theorem in \cite{11}, there is a unitary operator
\[
W: {\mathcal H}_{\text{sing}}(P(T \otimes I_n) \mid P{\mathcal H}^n + K') \oplus {\mathcal H} \to {\mathcal H}
\]
which intertwines $S_1 \oplus T$ and $T \mod {\mathcal C}_1$. The same kind of construction applied to $Q$ instead of $P$ produces a unitary operator
\[
V: {\mathcal H}_{\text{sing}}(Q(T \otimes I_n) \mid Q{\mathcal H}^n + K'') \oplus {\mathcal H} \to {\mathcal H}
\]
which intertwines $S_2 \oplus T$ and $T \mod {\mathcal C}_1$. On the other hand using the Kato--Rosenblum theorem $P(T \otimes I_n) \mid P{\mathcal H}^n + K'$ and $P(T \otimes I_n)P{\mathcal H}^n$ have unitarily equivalent absolutely continuous parts or equivalently the multiplicity function of these absolutely continuous parts are equal almost everywhere. Comparing this the corresponding situation for $Q$, we find there is a unitary operator
\[
U: {\mathcal H}_{ac}(P(T \otimes I_n) \mid P{\mathcal H}^n + K') \to {\mathcal H}_{ac}(Q(T \otimes I_n) \mid Q{\mathcal H}^n + K''
\]
which intertwines the restrictions of $P(T \otimes I_n) \mid P{\mathcal H}^n + K'$ and $Q(T \otimes I_n) \mid Q{\mathcal H}^n + K''$ to these subspaces. Then
\[
U \oplus V^*W: P{\mathcal H}^n \oplus {\mathcal H} \to Q{\mathcal H}^n \oplus {\mathcal H}
\]
is a unitary operator intertwining $(P(T \otimes I_n) \mid P{\mathcal H}^n + K') \oplus T$ and $(Q(T \otimes I_n) \mid Q{\mathcal H}^n + K'') \oplus T \mod {\mathcal C}_1$. But this then implies $U \oplus V^*W$ intertwines $P(T \otimes I_n) \mid P{\mathcal H}^n \oplus T$ and $Q(T \otimes I_n) \mid Q{\mathcal H}^n \oplus T \mod {\mathcal C}_1$. Hence by $2.6.$ we have $P \oplus I \sim Q \oplus I$.

Assume now that $P,Q \in {\mathcal P}_h({\mathcal E}(T \otimes I_n,{\mathcal C}_1)$ are so that $fmac(T)([P]_0) = fmac(T)([Q]_0)$. Then we have $mac(P(T \otimes I_n) \mid P{\mathcal H}^n) \mid {\mathbb R}\backslash\omega(T) = mac(Q(T \otimes I_n) \mid Q{\mathcal H}^n) \mid {\mathbb R}\backslash\omega(T)$. Adding $mac(T)\mid {\mathbb R}\backslash\omega(T)$ to both sides and then extending by $\infty$ over $\omega(T)$ then gives $mac({\tilde P}(T \otimes I_{n+1}) \mid {\tilde P}{\mathcal H}^n) = mac({\tilde Q}(T \otimes I_{n+1}) \mid {\tilde Q}{\mathcal H}^n)$ where ${\tilde P} = P \oplus I$, ${\tilde Q} = Q \oplus I$. By the first observation in this proof, this gives ${\tilde P} \oplus I \sim {\tilde Q} \oplus I$ which implies $[P]_0 = [Q]_0$.

The last result we proved actually implies $fmac(T)$ is injective, since $fmac(T)([P_1]_0 - [P_2]_0) = 0$ implies $fmac(T)([P_1]_0) = fmac(T)([P_2]_0)$ which then gives $[P_1]_0 = [P_2]_0$ and hence $[P_1]_0 - [P_2]_0 = 0$.\qed

\bigskip
\noindent
{\bf Remark 3.1.} The definition of $fmac(T)$ does not depend on the choice of $T$, that is, if we replace $T$ by $T'$ so that $T = T' \in {\mathcal C}_1$, then as remarked in $2.1.$ we have ${\mathcal E}(T,{\mathcal C}_1) = {\mathcal E}(T',{\mathcal C}_1)$ and the Kato--Rosenblum theorem immediately gives $fmac(T) = fmac(T')$. Remark also that ${\mathcal F}(T)$ identifies naturally with $K_0((T \mid {\mathcal H}_{ac}(T))')$, the prime denoting the commutant. Also the homomorphisms of $K_0$ we considered are ordered $K$-theory morphisms.

\bigskip
\noindent
{\bf Theorem 3.3.} {\em If $T$ is a Hermitian operator with $\sigma(T) = [0,1]$ and $({\mathcal J},|\ |_{\mathcal J})$ is a normed ideal so that ${\mathcal J} \ne {\mathcal C}_1$ and in which the finite rank operators are dense, then we have $K_0({\mathcal E}(T);{\mathcal J}) = 0$.
}

\bigskip
\noindent
{\bf {\em Proof.}} Since ${\mathcal J} \ne {\mathcal C}_1$, we have $k_{\mathcal J}(S) = 0$ for any Hermitian operator $S$, by what is essentially Kuroda's theorem (or use the more general \cite{1}). Remark also that if $P \in {\mathcal P}_h({\mathcal E}(T \otimes I_n;{\mathcal J}))$ then $\sigma(P(T \otimes I_n) \mid P{\mathcal H}^n) \subset [0,1]$. Thus we may use the adapted non-commutative Weyl--von~Neumann type theorem and get the existence of a unitary operator $V: P{\mathcal H}^n \oplus {\mathcal H} \to {\mathcal H}$ so that $((P(T \otimes I_n) \mid P{\mathcal H}^n) \oplus T)V - VT \in {\mathcal J}$. This implies $P \oplus I \sim (O \otimes I_n) \oplus I$ in ${\mathcal E}((T \otimes I_{n+1});{\mathcal M}_{n+1}{\mathcal J})$ and hence $[P]_0 + [I]_0 = [I]_0$ which gives $[P]_0 = 0$.\qed

\bigskip
\noindent
{\bf Remark 3.3.} The result in Theorem~$3.3$ is likely to hold more generally for $T$ a Hermitian operator with $\sigma(T) = \sigma_e(T)$. The proof would require a result for the ideal ${\mathcal J}$ similar to the result in \cite{4} which we used in the proof of Theorem~$3.2$, that is: if $T = T^*$, $\sigma(T) = \sigma_e(T)$ and $K = K^* \in {\mathcal J}$ then if $\lambda_j$, $j \in {\mathbb N}$ are the finite-multiplicity isolated eigenvalues of $T+K$ repeated according to multiplicity, then the diagonal operator with entries $d(\lambda_j,\sigma(T))$, $j \in {\mathbb N}$ is in ${\mathcal J}$. Such a result is probably known.

\section{Results in the case of an $n$-tuple of commuting Hermitian operators}
\label{sec4}

To extend the results of the preceding section to $n$-tuples of commuting Hermitian operators $\tau$, there is a new difficulty: the compression $P(\tau \otimes I_n) \mid P{\mathcal H}^n$ is no longer an $n$-tuple of commuting operators. Using a different approach it will be still possible to define a homomorphism $fmac(\tau)$ generalizing $fmac(T)$ of Theorem~$3.1$. We will also have to require that $n \ge 3$, because some of our results in \cite{11} which will be used here, have been established only for $n \ge 3$ in the strong version we need.

Let $\tau$ be an $n$-tuple of commuting Hermitian operators with $n \ge 3$ and assume that ${\mathcal H}_{ac}(\tau) = {\mathcal H}$, that is that the spectral measure of $\tau$ is Lebesgue absolutely continuous. By ${\mathcal F}(\tau)$ we denote the additive group of equivalence classes of almost everywhere equal maps $f: {\mathbb R}^n\backslash\omega(\tau) \to {\mathbb Z}$, so that $|f(x)| \le C\ mac(\tau)$ for some $C > 0$. Let also $i(\tau): (\tau)' \to {\mathcal E}(\tau;{\mathcal C}_n^-)$ be the inclusion of the commutant of $\tau$ into ${\mathcal E}(\tau;{\mathcal C}_n^-)$. From the fact that equivalence of projections in ${\mathcal P}_h({\mathcal M}_n((\tau)'))$ is the same as the equality of their central traces, one infers that there is a unique isomorphism $\alpha(\tau): K_0((\tau)') \to {\mathcal F}(\tau)$ so that $\alpha(\tau)([P]_0) = mac(P(\tau \otimes I_n) \mid P{\mathcal H}^n) \mid {\mathbb R}^n\backslash\omega(\tau) \in {\mathcal F}(\tau)$, when $P \in {\mathcal P}_h({\mathcal M}_n((\tau)'))$.

Let $u_m = f_m(\tau)$, $m \in {\mathbb N}$ where $f_m \in C^{\infty}({\mathbb R}^n)$ be unitary operators so that $u_m \to 0$ in the weak operator topology as $m \to \infty$. By Corollary~$1.6$ of (\cite{11} II) the strong limit $s - \lim_{m \to \infty} u^*_mXu_m = \psi(\tau)(X)$ exists for all $X \in {\mathcal E}(\tau;{\mathcal C}^-_n)$ and does not depend on the choice of the sequence $u_m$, $m \in {\mathbb N}$. In particular, replacing $u_m$ by $uu_m$, where $u = f(\tau)$ for some $f \in C^{\infty}({\mathbb R}^n)$ is unitary, we find that $u\psi(\tau)(X)u^* = \psi(\tau)(X)$ and hence $\psi(\tau)(X) \in (\tau)'$ for all $X \in {\mathcal E}(\tau;{\mathcal C}^-_n)$. Clearly if $X \in (\tau)'$ then we have $\psi(\tau)(X) = X$. This shows that $\psi(\tau) \circ i(\tau) = id_{(\tau)'}$, so that $(\tau)'$ is a retract of ${\mathcal E}(\tau;{\mathcal C}^-_n)$.

Using $\alpha(\tau)$ to identify ${\mathcal F}(\tau)$ and $K_0((\tau)')$ we get that the homomorphism
\[
K_0(i(\tau)) \circ (\alpha(\tau))^{-1}: {\mathcal F}(\tau) \to K_0({\mathcal E}(\tau;{\mathcal C}_n^-))
\]
has a right inverse $fmac(\tau)$ which equals
\[
\alpha(\tau) \circ K_0(\psi(\tau)): K_0({\mathcal E}(\tau;{\mathcal C}^-_n)) \to {\mathcal F}(\tau)
\]
and which makes ${\mathcal F}(\tau)$ a retract of $K_0({\mathcal E}(\tau;{\mathcal C}^-_n))$.

Remark also that if ${\tilde \tau}$ is another $n$-tuple of commuting Hermitian operators so that ${\tilde T}_j - T_j \in {\mathcal C}^-_n$, $1 \le j \le n$ and ${\mathcal H}_{ac}({\tilde \tau}) = {\mathcal H}$, then performing all the above constructions with ${\tilde \tau}$ instead of $\tau$ we get the same homomorphism $fmac$, that is $fmac({\tilde \tau}) = fmac(\tau)$. Indeed, notice first that ${\mathcal E}(\tau;{\mathcal C}^-_n) = {\mathcal E}({\tilde \tau};{\mathcal C}^-_n)$ and that by (\cite{11} II) there is a unitary operator $W$ so that $W \tau W^* = {\tilde \tau}$ and that this implies $W \in {\mathcal E}(\tau;{\mathcal C}^-_n)$. This gives $W(\tau)'W^* = ({\tilde \tau})'$ and $\psi({\tilde \tau})(X) = W(\psi(\tau)(W^*XW))W^*$. Moreover we have $\alpha({\tilde \tau})([(W \otimes I_n)(P)(W^* \otimes I_n)]_0) = \alpha(\tau)([P]_0)$ if $P \in {\mathcal P}_h({\mathcal M}_n((\tau)')$. We then see that
\[
\begin{aligned}
K_0(i({\tilde \tau})) \circ \alpha({\tilde \tau})^{-1} &= K_0(\text{Ad } W) \circ K_0(i(\tau)) \circ \alpha(\tau)^{-1} \\
&= K_0(i(\tau)) \circ \alpha(\tau)^{-1}
\end{aligned}
\]
and
\[
\begin{aligned}
fmac(\tau) &= \alpha(\tau) \circ K_0(\psi(\tau)) \\
&= \alpha({\tilde \tau}) \circ K_0(\text{Ad } W \circ \psi(\tau)) \\
&= \alpha({\tilde \tau}) \circ K_0(\psi({\tilde \tau}) \circ \text{Ad } W) \\
&= \alpha({\tilde \tau}) \circ K_0(\psi({\tilde \tau}) \circ K_0(\text{Ad } W) \\
&= \alpha({\tilde \tau}) \circ K_0(\psi({\tilde \tau}) = fmac({\tilde \tau}).
\end{aligned}
\]

We should also remark that $\alpha(\tau)$ and $K_0(\psi(\tau))$ are also compatible with the order structure on the $K_0$ groups.

We shall summarize the above discussion as the next theorem.

\bigskip
\noindent
{\bf Theorem 4.1.} {\em Let $\tau$ be an $n$-tuple of commuting Hermitian operators, $n \ge 3$ and satisfying ${\mathcal H}_{ac}(\tau) = {\mathcal H}$. Then there is an isomorphism $\alpha(\tau): K_0((\tau)') \to {\mathcal F}(\tau)$ and a $*$-homomorphism $\psi(\tau): {\mathcal E}(\tau;{\mathcal C}^-_n) \to (\tau)'$ defined by $\psi(X) = s - \lim_{m \to \infty} u_mXu*_m$ where $u_m = f_m(\tau)$, $f_m \in C^{\infty}({\mathbb R}^n)$ are unitary operators converging weakly to $O$. With $i(\tau)$ denoting the inclusion $(\tau)' \hookrightarrow {\mathcal E}(\tau;{\mathcal C}^-_n)$, the homomorphism $fmac(\tau) = \alpha(\tau) \circ K_0(\psi(\tau))$ of $K_0({\mathcal E}(\tau;{\mathcal C}^-_n))$ to ${\mathcal F}(\tau)$, is a retraction for the inclusion ${\mathcal F}(\tau) \hookrightarrow K_0({\mathcal E}(\tau;{\mathcal C}^-_n))$, thus making ${\mathcal F}(\tau)$ a retract of $K_0({\mathcal E}(\tau;{\mathcal C}^-_n))$. Replacing $\tau$ with ${\tilde \tau}$ so that $\tau \equiv {\tilde \tau} \mod {\mathcal C}^-_n$ and ${\mathcal H}_{ac}({\tilde \tau}) = {\mathcal H}$, yields the same homomorphism, that is $fmac(\tau) = fmac({\tilde \tau})$ and ${\mathcal F}(\tau) = {\mathcal F}({\tilde \tau})$ identifies in the same way with the same direct summand of $K_0({\mathcal E}(\tau;{\mathcal C}^-_n))$.
}

\bigskip
\noindent
{\bf Remark 4.1.} The case $n = 2$ is a good occasion to reflect about how incomplete our understanding of $K_0({\mathcal E}(\tau;{\mathcal C}^-_n))$ actually is if $n > 1$. On one hand the construction of $\psi(\tau)$ is not available since in \cite{11} II we obtained the strong results about generalized wave operators only for $n \ge 3$ and we don't know whether the weaker results we proved for $n = 2$ can be strengthened. On the other hand there are homomorphisms related to cyclic homology like the one in \cite{14} for $K_0({\mathcal E}(\tau;{\mathcal C}_2))$ where $\tau = (T_1,T_2)$ which associates to $[P]_0$ the Pincus $g$-function of $P(T_1 + iT_2)P$ and which can be composed with the $K_0$-map for the inclusion ${\mathcal E}(\tau;{\mathcal C}^-_2) \hookrightarrow {\mathcal E}(\tau;{\mathcal C}_2)$.

\section{The ideal of compact operators}
\label{sec5}

Assuming $k_{\mathcal I}(\tau) = 0$ in \cite{14}, \cite{15} we showed that the ideal of compact operators ${\mathcal K}(\tau;{\mathcal I})$ and the quotient $\xi/{\mathcal K}(\tau;{\mathcal I})$ have several nice properties, in particular ${\mathcal E}/{\mathcal K}(\tau;{\mathcal I})$ is isometric to its image in the Calkin algebra which is a $C^*$-algebra. Since many of the $K$-theory results in this paper are about ${\mathcal E}(\tau;{\mathcal I})$ with $k_{\mathcal I}(\tau) > 0$, we explore in this section what happens with ${\mathcal K}(\tau;{\mathcal I})$ and ${\mathcal E}/{\mathcal K}(\tau;{\mathcal I})$ if $k_{\mathcal I}(\tau) > 0$. It turns out that if $k_{\mathcal I}(\tau) < \infty$ there are weaker forms of some of the results which still hold. It may be interesting to recall that when ${\mathcal I} = {\mathcal C}_{\infty}^-$ we always have $k_{\infty}^-(\tau) < \infty$.

\bigskip
\noindent
{\bf Proposition 5.1.} {\em Assume that ${\mathcal R}$ is dense in ${\mathcal I}$ and that $k_{\mathcal I}(\tau) = C < \infty$. Let $A_n \in {\mathcal R}_1^+$ be so that $A_n \uparrow I$ as $n \to \infty$ and that $\lim_{n \to \infty}|[A_n,\tau]|_{\mathcal I} = C$. Then we have $\lim_{n \to \infty} |\|A_n\|| = 2 + C$ and
\[
\lim_{n \to \infty} |\|K - A_nK\|| = \lim_{n \to \infty} |\|K - KA_n\|| = 0
\]
when $K \in {\mathcal K}(\tau;{\mathcal I})$.
}

\bigskip
\noindent
{\bf {\em Proof.}} We have $\lim_{n \to \infty} \|A_n\| = 1$ and $\lim_{n \to \infty} |[A_n,\tau]|_{\mathcal I} = C$ so that indeed $\lim_{n \to \infty}|\||A_n\|| = 1 + C$.

If $K \in {\mathcal K}(\tau;{\mathcal I})$ using only the fact that $K$ is a compact operator we have
\[
\lim_{n \to \infty} \|(I-A_n)K\| = 0.
\]

On the other hand
\[
|[(I-A_n)K,\tau]|_{\mathcal I} \le |[A_n,\tau]|_{\mathcal I}\|K\| + |(I-A_n)[K,\tau]|_{\mathcal I}.
\]

Since $[K,T_j] \in {\mathcal I}$ and ${\mathcal R}$ is dense in ${\mathcal I}$ we have that
\[
\lim_{n \to \infty} |(I-A_n)[K,\tau]|_{\mathcal I} = 0
\]
and hence
\[
\limsup_{n \to \infty}\||(I-A_n)K|\| \le C\|K\|.
\]
This last inequality can be improved. Indeed, if $R \in {\mathcal R}$ we have
\[
|[(I-A_n)R,\tau]|_{\mathcal I} \le 2(\mbox{rank } R) \cdot \|[(I-A_n)R,\tau]\|
\]
so that
\[
\lim_{n \to \infty} |\|(I-A_n)R\|| = 0
\]
which then gives
\[
\limsup_{n \to \infty} |\|(I-A_n)K\|| \le \limsup_{n \to \infty} |\|(I-A_n)(K-R)\|| \le C \|-R\|.
\]

Since $R \in R$ is arbitrary we get that
\[
\lim_{n \to \infty} |\|(I-A_n)K\|| = 0.
\]
Replacing $K$ by $K^*$ we also get
\[
\lim_{n \to \infty} |\|K(I-A_n)\|| = 0.
\]
\qed

\bigskip
\noindent
{\bf Corollary 5.1.} {\em If $k_{\mathcal I}(\tau) < \infty$ and $R$ is dense in ${\mathcal I}$, then $R$ is dense in ${\mathcal K}(\tau;{\mathcal I})$.
}

\bigskip
\noindent
{\bf Corollary 5.2.} {\em If $k_{\mathcal I}(\tau) < \infty$ and $R$ is dense in ${\mathcal I}$, then $K_0({\mathcal K}(\tau;{\mathcal I})) = {\mathbb Z} \cdot [P]_0$ where $P$ is a rank one projection and $K_1({\mathcal K}(\tau;{\mathcal I})) = 0$.
}

\bigskip
\noindent
{\bf Corollary 5.3.} {\em If $k_{\mathcal I}(\tau) < \infty$ and $R$ is dense in ${\mathcal I}$, then the algebraic identification of ${\mathcal E}/{\mathcal K}(\tau;{\mathcal I})$ and $p({\mathcal E}(\tau;{\mathcal I})) \subset B/{\mathcal K}$ is an isomorphism of Banach $*$-algebras and $p({\mathcal E}(\tau;{\mathcal I}))$ is a $C^*$-algebra. Moreover we have
\[
\|p(T)\| \le \|T + {\mathcal K}(\tau;{\mathcal I})\|_{{\mathcal E}/{\mathcal K}} \le (1 + k_{\mathcal I}(\tau))\|p(T)\|
\]
if $T \in {\mathcal E}(\tau;{\mathcal I})$.
}

\bigskip
\noindent
{\bf {\em Proof.}} If $T \in {\mathcal E}(\tau;{\mathcal I})$, with the notation of Proposition~5.1 we have $\|(I-A_n)T\| \to \|p(T)\|$ as $n \to \infty$. We also have
\[
|[(I-A_n)T,\tau]|_{\mathcal I} \le |[A_n,\tau]|_{\mathcal I}\|T\| + |(I-A_n)[T,\tau]|_{\mathcal I}.
\]
Since $R$ is dense in ${\mathcal I}$ we have that
\[
\lim_{n \to \infty}|(I-A_n)[T,\tau]|_{\mathcal I} = 0.
\]
Thus we get
\[
\limsup_{n \to \infty}|[(I-A_n)T,\tau]|_{\mathcal I} \le k_{\mathcal I}(\tau)\|T\|
\]
which in view of
\[
\limsup_{n \to \infty} |[(I-A_n)R,\tau]|_{\mathcal I} = 0
\]
if $R \in R$, gives
\[
\limsup_{n \to \infty} |[(I-A_n)T,\tau]|_{\mathcal I} \le k_{\mathcal I}(\tau)\|T-R\|
\]
and $R \in R$ being arbitrary we get
\[
\limsup_{n \to \infty} |[(I-A_n)T,\tau]|_{\mathcal I} \le k_{\mathcal I}(\tau)\|p(T)\|.
\]
This implies
\[
\|T+{\mathcal K}(\tau;{\mathcal I})\|_{{\mathcal E}/{\mathcal K}} \le \limsup_{n \to \infty} |\|(I-A_n)\tau\|| \le (1 + k_{\mathcal I}(\tau))\|p(T)\|.
\]
On the other hand clearly
\[
\|T+{\mathcal K}(\tau;{\mathcal I}\|_{{\mathcal E}/{\mathcal K}} = \inf_{K \in {\mathcal K}(\tau;{\mathcal I})} |\|T-K\|| \ge \inf_{K \in {\mathcal K}(\tau;{\mathcal I})} \|T-K\| \ge \|p(T)\|.
\]
The equivalence of $\|p(T)\|$ and $\|T+{\mathcal K}(\tau;{\mathcal I})\|_{{\mathcal E}/{\mathcal K}}$ then easily gives the remaining assertions of the Corollary. \qed

\bigskip
\noindent
{\bf Remark 5.1.} The Macaev ideal ${\mathcal C}_{\infty}^-$ has the remarkable property that $k_{\infty}^-(\tau) < \infty$ for all $n$-tuples $\tau$ and moreover $R$ is dense in ${\mathcal C}_{\infty}^-$. Thus if ${\mathcal I} = {\mathcal C}_{\infty}^-$ the assumption of Proposition~5.1 and its corollaries Corollary~5.1, Corollary~5.2, Corollary~5.3 is satisfied for $\tau$. We record part of this in the next corollary.

\bigskip
\noindent
{\bf Corollary 5.4.} {\em If $\tau$ is a $n$-tuple of operators then $R$ is dense in ${\mathcal K}(\tau;{\mathcal C}_{\infty}^-)$ so that $K_0({\mathcal K}(\tau;{\mathcal C}_{\infty}^-)) \simeq {\mathbb Z}$ and $K_1({\mathcal K}(\tau;{\mathcal C}_{\infty}^-)) = 0$. Moreover $p({\mathcal E}(\tau;{\mathcal C}_{\infty}^-))$ is a $C^*$-subalgebra of the Calkin algebra $B/{\mathcal K}$ which is isomorphic as a Banach $*$-algebra to ${\mathcal E}/{\mathcal K}(\tau;{\mathcal C}_{\infty}^-)$.
}

\bigskip
\noindent
{\bf Corollary 5.5.} {\em Assume that $k_{\mathcal I}(\tau) < \infty$, that $R$ is dense in ${\mathcal I}$ and that the $C^*$-algebra $C^*(\tau)$ generated by $\tau$ and $I$ has a one-dimensional representation. Then $p({\mathcal E}(\tau;{\mathcal I}))$ is a $C^*$-algebra which contains a Fredholm element of index~$1$ and we have canonical isomorphisms $K_0({\mathcal E}(\tau;{\mathcal I})) \simeq K_0({\mathcal E}/{\mathcal K}(\tau;{\mathcal I}))$ and $K_1({\mathcal E}(\tau;{\mathcal I})) \simeq \ker \partial$ where $\partial$ is the index map $K_1({\mathcal E}/{\mathcal K}(\tau;{\mathcal I})) \to K_0({\mathcal K}(\tau;{\mathcal I})) \simeq {\mathbb Z}$.
}

\bigskip
\noindent
{\bf {\em Proof.}} The assumptions imply that $\tau$ is unitarily equivalent $\mod {\mathcal I}$ (actually $\mod {\mathcal C}_1$) to $\tau \oplus (\lambda_jI)_{1 \le j \le n}$ for some complex numbers $\lambda_j$. The commutant of $\tau \oplus (\lambda_jI)_{1 \le j \le n}$ clearly contains a Fredholm operator of index~$1$. This implies that the index map $K_1({\mathcal E}/{\mathcal K}(\tau;{\mathcal I})) \to K_0({\mathcal K}(\tau;{\mathcal I})) \simeq {\mathbb Z}$ is surjective and since we also have $K_1({\mathcal K}(\tau;{\mathcal I})) = 0$ we get along standard lines the assertions.\qed

\bigskip
Since Corollary~5.1 gives conditions under which $R$ is dense in ${\mathcal K}(\tau;{\mathcal I})$, we will conclude this section pointing a few results we have established previously under assumptions which included $k_{\mathcal I}(\tau) = 0$ and which can now be seen to extend to situations when $k_{\mathcal I}(\tau) < \infty$.

\bigskip
\noindent
{\bf Proposition 5.2.} {\em Assume $k_{\mathcal I}(\tau) < \infty$ and $R$ is dense in ${\mathcal I}$. Then the dual of ${\mathcal K}(\tau;{\mathcal I})$ identifies isometriclly with $({\mathcal C}_1 \times ({\mathcal I}^*)^n)/{\mathcal N}$ where the norm on ${\mathcal C}_1 \times ({\mathcal I}^*)^n$ is
\[
\|(x,(y_j)_{1 \le j \le n})\| = \max\left(|x|_1,\sum_{1 \le j \le n} |y_j|_{{\mathcal I}^*}\right)
\]
and
\[
{\mathcal N} = \left\{\left( \sum_{1 \le j \le n} [T_j,y_j],(y_j)_{1 \le j \le n}\right) \in {\mathcal C}_1 \times ({\mathcal I}^*)^n \mid \sum_{1 \le j \le n} [T_j,y_j] \in {\mathcal C}_1\right\},
\]
the duality pairing arising from
\[
(K,(x,(y_j)_{1 \le j \le n})) \to \mbox{\em Tr}\left(Kx + \sum_{1 \le j \le n} [T_j,K]_{y_j}\right).
\]
}

\bigskip
The proof of the preceding proposition is the same as that of Proposition~4.3 in \cite{15} and uses the density of $R$ in ${\mathcal K}(\tau;{\mathcal I})$ to identify ${\mathcal N}$ as the annihilator of ${\mathcal K}(\tau;{\mathcal I})$ in ${\mathcal C}_1 \times ({\mathcal I}^*)^n$.

Let $\Phi$ be a norming function and ${\mathcal G}_{\Phi}$ the largest ideal determined by $\Phi$ and ${\mathcal G}_{\Phi}^{(0)}$ the closure of $R$ in ${\mathcal G}_{\Phi}$. Then Lemma~4.1 and Proposition~4.2 of \cite{15} are easily extended to the case when $k_{\Phi}(\tau) < \infty$. Thus ${\mathcal K}(\tau;{\mathcal G}_{\Phi}^{(0)})$ is a closed two-sided ideal in ${\mathcal E}(\tau;{\mathcal G}_{\Phi})$, the unit ball of ${\mathcal E}(\tau;{\mathcal G}_{\Phi})$ is weakly compact and ${\mathcal M}({\mathcal K}(\tau;{\mathcal G}_{\Phi}^{(0)}))$ the algebra of bounded multipliers of ${\mathcal K}(\tau;{\mathcal G}_{\Phi}^{(0)})$ identifies with ${\mathcal E}(\tau;{\mathcal G}_{\Phi})$.

\end{document}